\documentclass[a4paper,11pt,reqno]{amsart}

\usepackage{latexsym,dsfont,amssymb,amsmath,amsthm,amscd}

\usepackage[all]{xy}

\chardef\bslash=`\\ 

\numberwithin{equation}{section}

\newtheorem{theorem}{Theorem}[section]
\newtheorem{corollary}[theorem]{Corollary}
\newtheorem{lemma}[theorem]{Lemma}

\theoremstyle{remark}
\newtheorem{remark}[theorem]{Remark}

\theoremstyle{definition}

\newcommand\bp{\begin{proof}}
\newcommand\ep{\end{proof}}

\newcommand\Dhat{{\hat\Delta}}

\newcommand\A{{\mathcal A}}
\newcommand\CC{{\mathcal C}}
\newcommand\D{{\mathcal D}}
\newcommand\E{{\mathcal E}}
\newcommand\F{{\mathcal F}}

\newcommand\RR{{\mathcal R}}
\newcommand\U{{\mathcal U}}

\newcommand\g{{\mathfrak g}}
\newcommand\h{{\mathfrak h}}

\newcommand\sltwo{\mathfrak{sl}_2}

\newcommand\Clg{\operatorname{Cl}({\mathfrak g})}

\newcommand\Dom{\operatorname{Dom}}

\newcommand\End{\operatorname{End}}

\newcommand\Tr{\operatorname{Tr}}

\newcommand\Vect{\mathcal Vec}

\newcommand{\ad}{\operatorname{ad}}
\newcommand{\add}{{\widetilde\ad}}

\newcommand{\C}{{\mathbb C}}

\newcommand{\N}{{\mathbb N}}

\newcommand{\R}{{\mathbb R}}
\newcommand\Sp{{\mathbb S}}
\newcommand\T{{\mathbb T}}
\newcommand\Z{{\mathbb Z}}

\newcommand\eps{\varepsilon}

\newcommand\enu[1]{\smallskip\newline\makebox[5mm][l]{\rm(#1)}}

\begin{document}

\title[K-homology class of the Dirac operator]
{K-homology class of the Dirac operator on a compact quantum group}


\author[S. Neshveyev]{Sergey Neshveyev}
\address{Department of Mathematics, University of Oslo,
P.O. Box 1053 Blindern, NO-0316 Oslo, Norway}

\email{sergeyn@math.uio.no}

\author[L. Tuset]{Lars Tuset}
\address{Faculty of Engineering, Oslo University College,
P.O. Box 4 St.~Olavs plass, NO-0130 Oslo, Norway}
\email{Lars.Tuset@iu.hio.no}

\thanks{Supported by the Research Council of Norway.}

\date{February 1, 2011}

\begin{abstract}
By a result of Nagy, the C$^*$-algebra of continuous functions on the $q$-deformation $G_q$ of a simply connected semisimple compact Lie group $G$ is KK-equivalent to $C(G)$. We show that under this equivalence the K-homology class of the Dirac operator on $G_q$, which we constructed in an earlier paper, corresponds to that of the classical Dirac operator. Along the way we prove that for an appropriate choice of isomorphisms between completions of~$U_q\g$ and $U\g$ a family of Drinfeld twists relating the deformed and classical coproducts can be chosen to be continuous in $q$.
\end{abstract}

\maketitle

\bigskip

\section*{Introduction}

In \cite{NT2} we constructed a Dirac operator $D_q$ on the $q$-deformation $G_q$ of any simply connected semisimple compact Lie group $G$. The construction involved a special unitary element $\F^q$ in the von Neumann algebra $W^*(G)\bar\otimes W^*(G)$, which relates the coproducts in $W^*(G_q)$ and $W^*(G)$. The existence of such an element, called a unitary Drinfeld twist, is a consequence of a fundamental and highly nontrivial result in quantum group theory due to Kazhdan and Lusztig~\cite{KL1}, see also \cite{NT3}. Since the construction of a Drinfeld twist is involved and not particularly explicit, certain properties of the operators $D_q$ are not immediate. In particular, even though it is intuitively clear that $D_q$ is a deformation of the classical Dirac operator and therefore should in some sense define the same index map on K-theory, it is not even obvious that the K-homology class of~$D_q$ is always nonzero. The goal of this note is to show that it is indeed nonzero. In fact we show that the K-homology class of $D_q$
corresponds exactly to that of the classical Dirac operator via the KK-equivalence of~$C(G_q)$ and~$C(G)$ established by Nagy~\cite{Na}.
Therefore, upon identifying the K-theories of $G_q$ and $G$, the index map defined by $D_q$ does not depend on $q$, as expected. We should remark that the question of invariance of the index map under deformation has been studied in a recent paper by Yamashita~\cite{Ya} in the context of Connes-Landi $\theta$-deformations.

\smallskip

The paper is organized as follows. In Section~\ref{sconfield} we show that the family of C$^*$-algebras $C(G_q)$ has a canonical continuous field structure. The result is more or less known~\cite{Na}, but is usually formulated in terms of standard generators of $\C[G_q]$. We propose a simpler approach based on a natural notion of a continuous family of isomorphisms $W^*(G_q)\cong W^*(G)$.

In Section~\ref{stwist} we prove that once a continuous family of isomorphisms $W^*(G_q)\cong W^*(G)$ is fixed, the corresponding family of Drinfeld twists $\F^q$ can be chosen to be continuous in $q$. This result is not strictly speaking necessary for our main result on $D_q$, for which it suffices to know that $D_q$ does not depend on $\F^q$ for a fixed isomorphism $W^*(G_q)\cong W^*(G)$, see~\cite{NT4}, but it simplifies the arguments and is of independent interest. Both results, continuity of $\F^q$ and uniqueness of $D_q$, depend crucially on the fact that any two unitary Drinfeld twists differ by the coboundary of a central unitary element, a result we proved in~\cite{NT4}.

In Section~\ref{sdirac} we prove our main result. For this we show that the family of operators $D_q$ define a Kasparov module for the algebra of continuous sections of $(C(G_q))_{q\in[a,b]}$ and the algebra $C[a,b]$. With the preparation in the previous two sections the proof essentially boils down to observing that estimates in our paper~\cite{NT2} are uniform in ~$q$.

\medskip\noindent
{\bf Acknowledgement.} The authors are grateful to Christian Voigt and Makoto Yamashita for stimulating discussions.


\bigskip

\section{Continuous field of function algebras}
\label{sconfield}

Let $G$ be a simply connected semisimple compact Lie group, $\g$ its
complexified Lie algebra, $\h\subset\g$ the Cartan subalgebra defined by a maximal torus in $G$. Fix a system $\{\alpha_1,\dots,\alpha_r\}$ of simple roots. Let $(a_{ij})_{1\le i,j\le r}$ be the Cartan matrix defined by $\alpha_1,\dots,\alpha_r$, and $d_1,\dots,d_r$ be
the coprime positive integers such that $(d_ia_{ij})_{i,j}$ is
symmetric. 

For every $q>0$, $q\ne 1$, consider the quantized universal enveloping algebra $U_q\g$ with generators $E_i^q,F_i^q,K_i^q$ ($1\le i\le r$); we follow the conventions in~\cite{NT3,NT4}. For $q=1$ we let $U_1\g=U\g$ and denote by $E^1_i,F^1_i,h^1_i$ the standard generators of $U\g$. We will often omit various indices corresponding to~$q=1$. Consider the category $\CC_q(\g)$ of finite dimensional admissible $U_q\g$-modules and denote by $\U(G_q)$ the endomorphism ring of the forgetful functor $\CC_q(\g)\to\Vect$. We think of $\U(G_q)$ as a completion of $U_q\g$. For $q\ne1$ denote by $h^q_i$ the unique self-adjoint element of $\U(G_q)$ such that $K^q_i=q^{d_i h^q_i}$.

Since simple objects of $\CC_q(\g)$ are classified by dominant integral weights, for every $q>0$ we have canonical identifications of the centers of $\U(G_q)$ and $\U(G)$. It extends to a $*$-isomorphism $\varphi^q\colon \U(G_q)\to \U(G)$. For every dominant integral weight $\lambda\in P_+$ and $q>0$ fix an irreducible $*$-representa\-tion $\pi^q_\lambda\colon U_q\g\to B(V^q_\lambda)$ with highest weight $\lambda$ and a highest weight unit vector $\xi^q_\lambda\in V^q_\lambda$. Then to define $\varphi^q$ is the same as to fix an isomorphism $B(V^q_\lambda)\cong B(V_\lambda)$ for every $\lambda$.

\smallskip

We say that a family $\{\varphi^q\}_{q>0}$ of $*$-isomorphisms extending the canonical identifications of the centers is continuous if for every finite dimensional representation $\pi$ of $U\g$ the operators $\pi(\varphi^q(E^q_i))$, $\pi(\varphi^q(F^q_i))$ and $\pi(\varphi^q(h^q_i))$ depend continuously on~$q$.

\begin{lemma}
A continuous family of $*$-isomorphisms $\varphi^q\colon \U(G_q)\to \U(G)$ always exists.
\end{lemma}

\bp  It suffices to show that for every $\lambda\in P_+$ there exist unitaries $u^q_\lambda \colon V^q_\lambda\to V_\lambda$ such that the operators $u^q_\lambda\pi^q_\lambda(X^q){u^q_\lambda}^*$ depend continuously on $q$ for $X^q=E^q_i,F^q_i,h^q_i$. For this, in turn, it is enough to show that such unitaries exist locally.

Therefore fix $\lambda$ and $q_0>0$. For every multi-index $I=(i_1,\dots,i_k)$ ($1\le i_j\le r$) define
$$
e^q_I=\pi^q_\lambda(F^q_{i_1}\dots F^q_{i_k})\xi^q_\lambda\in V^q_\lambda.
$$
We can choose multi-indices $I_1,\dots, I_n$ such that the vectors $e^{q_0}_{I_1},\dots,e^{q_0}_{I_n}$ form a basis in~$V^{q_0}_\lambda$. The quantum Serre relations, together with the identities ${F_i^q}^*=(K_i^q)^{-1}E_i^q$, imply that the scalar products~$(e^q_I,e^q_J)$ depend continuously on $q$. Hence, for some $\eps>0$, the vectors $e^{q}_{I_1},\dots,e^{q}_{I_n}$ form a basis in~$V^{q}_\lambda$ for all $q\in(q_0-\eps,q_0+\eps)$. Applying the Gram-Schmidt orthogonalization we get an orthonormal basis $\zeta^q_1,\dots,\zeta^q_n$ in~$V^{q}_\lambda$. Let $v^q\colon V^q_\lambda\to V^{q_0}_\lambda$ be the unitary mapping $\zeta^q_i$ into $\zeta^{q_0}_i$. By construction, for every multi-index~$I$, the coefficients of $e^q_I$ in the basis $\zeta^q_1,\dots,\zeta^q_n$ depend continuously on $q$, hence the matrix coefficients of $\pi^q_\lambda(F_i^q)$ in this basis also depend continuously on $q$. It follows that the operators $v^q\pi^q_\lambda(F_i^q){v^q}^*$ depend continuously on $q$. The same is clearly true for $h^q_i$  in place of $F^q_i$ (in fact, we even have $v^q\pi^q_\lambda(h_i^q){v^q}^*=h_i^{q_0}$), hence also for $E^q_i=q^{d_ih^q_i}{F_i^q}^*$.

Now take an arbitrary unitary $u\colon V^{q_0}_\lambda\to V_\lambda$. Then the unitaries $u_\lambda^q=uv^q$, $q\in(q_0-\eps,q_0+\eps)$, have the required properties.
\ep

From now on we will fix a continuous family of $*$-isomorphisms $\varphi^q\colon \U(G_q)\to \U(G)$ such that $\varphi^1=\iota$.

\smallskip

For every $q>0$ denote by $\C[G_q]\subset\U(G_q)^*$ the Hopf $*$-algebra of matrix coefficients of finite dimensional admissible $U_q\g$-modules, and by $C(G_q)$ its C$^*$-completion.

\begin{theorem} \label{tcontfield}
The family of C$^*$-algebras $C(G_q)$, $q>0$, has a unique structure of a continuous field of C$^*$-algebras such that for every $a\in\C[G]$ the section $q\mapsto a\varphi^q\in C(G_q)$ is continuous. This structure does not depend on the choice of a continuous family of $*$-isomorphisms~$\varphi^q$.
\end{theorem}

It is known that the C$^*$-algebras $C(SU_q(N))$ have a continuous field structure such that the matrix coefficients of the fundamental representation form continuous sections. For $N=2$ this was proved by Bauval~\cite{Bau} and for all $N\ge2$ by Nagy~\cite{Na}. As mentioned in~\cite{Na} the same proof as for~$SU_q(N)$ works for all other classical simple compact Lie groups. In principle the same result is also true for exceptional groups once explicit generators of $\C[G_q]$ have been found. The point of the above theorem is that there is actually no need to do this, it is enough to know that there exists a `continuous' choice of generators.

Although we will not need this here, we note that the theorem and its proof also imply that the families of function algebras on $q$-deformations of homogeneous spaces of $G$ can be given a continuous field structure without working out explicit generators and relations in those algebras. This, as well as the relation of the above result to Rieffel's notion of strict deformation quantization, will be discussed elsewhere.

\bp[Proof of Theorem~\ref{tcontfield}] First consider the dependence of the continuous field structure on the isomorphisms $\varphi^q$. Assume we have another continuous family of isomorphisms $\psi^q$. For every $\lambda\in P_+$ denote by $\gamma^q_\lambda$ the unique automorphism of $B(V_\lambda)$ such that $\pi_\lambda\psi^q=\gamma^q_\lambda\pi_\lambda\varphi^q$. Then the map $q\mapsto\gamma^q_\lambda\in\operatorname{Aut}(B(V_\lambda))$ is continuous. It follows that for any linear functional $\omega$ on $B(V_\lambda)$ the elements $\omega\gamma^q\pi_\lambda\in \C[G]$ decompose into finite linear combinations of elements $\nu\pi_\lambda$ with continuous coefficients, so that the sections $q\mapsto \omega\pi_\lambda\psi^q$ are finite linear combinations of sections $q\mapsto\nu\pi_\lambda\varphi^q$ with continuous coefficients. Therefore if the latter sections are continuous, the former are continuous as well.

\smallskip

Since $\C[G_q]$ is dense in $C(G_q)$, it is also clear that the  continuous field structure is unique if it exists.

\smallskip

To prove existence, first consider the case $G=SU(2)$. As usual identify $P$ with the half-integers. For every $s\in\frac{1}{2}\Z_+$ consider the orthonormal basis in~$V^q_s$ consisting of the vectors $\|(F^q)^k\xi^q_s\|^{-1}(F^q)^k\xi^q_s$, $k=0,\dots,2s$. Let $u^{s;q}_{ij}\in\C[SU_q(2)]$ be the matrix coefficients of $\pi^q_s$ in this basis. We also use these bases to construct the isomorphisms $\varphi^q$, so that $u^{s;q}_{ij}
=u^{s}_{ij}\varphi^q$. By~\cite{Bau}, see also~\cite{Bl}, there exists a continuous field structure on the C$^*$-algebras $C(SU_q(2))$ such that the sections $q\mapsto u^{1/2;q}_{ij}$ are continuous. To prove that the sections $q\mapsto u^{s;q}_{ij}$ are continuous in this structure we just have to show that $u^{s;q}_{ij}$ can be expressed as polynomials of the elements $u^{1/2;q}_{ij}$ with continuous coefficients. That this is indeed possible is easy to see using the unique embedding of $V^q_s$ into $(V^q_{1/2})^{\otimes 2s}$ mapping~$\xi^q_s$ into $(\xi^q_{1/2})^{\otimes 2s}$; in fact an explicit expression for these polynomials is known~\cite{VS}.

\smallskip

Turning to the general case, for every simple root $\alpha_i$ consider the $*$-homomorphism $\sigma^q_i\colon C(G_q)\to C(SU_{q^{d_i}}(2))$ which is dual to the embedding $\rho_i^q\colon \U(SU_{q^{d_i}}(2))\hookrightarrow \U(G_q)$ defined by $E^{q^{d_i}}\mapsto E_i^q$, $F^{q^{d_i}}\mapsto F_i^q$ and $h^{q^{d_i}}\mapsto h_i^q$. Let $w=s_{i_1}\dots s_{i_n}$ be the longest element in the Weyl group of $G$ written in reduced form. Put $A_q=C(SU_{q^{d_{i_1}}}(2))\otimes\dots\otimes C(SU_{q^{d_{i_n}}}(2))$. Since $C(SU_q(2))$ is a C$^*$-algebra of type I, the C$^*$-algebra of continuous sections of the field $(C(SU_q(2)))_{q>0}$ vanishing at infinity is of type I as well, hence exact. By~\cite[Theorem~4.6]{KW} it follows that the field $(A_q)_{q>0}$ has a continuous field structure such that the tensor product of continuous sections is a continuous section. Define a $*$-homomorphism
$$
\sigma^q\colon C(G_q)\to A_q\ \ \hbox{by}\ \ \sigma^q(a)=(\sigma^q_{i_1}\otimes\dots\otimes \sigma^q_{i_n})\Delta^{(n-1)}_q(a).
$$
It follows from the description of irreducible representations of $C(G_q)$, see e.g.~\cite[Theorem~6.2.7]{KS}, that $\sigma^q$ is injective for every~$q$. Therefore the field $(C(G_q))_{q>0}$ embeds into $(A_q)_{q>0}$, so to prove existence of the required continuous field structure on $(C(G_q))_{q>0}$ it suffices to show that for every $a\in \C[G]$ the section $q\mapsto\sigma^q(a\varphi^q)\in A_q$ is continuous.

Since $\varphi^q$ is an algebra homomorphism, the dual map $\C[G]\to\C[G_q]$, $a\mapsto a\varphi^q$, is a coalgebra homomorphism. Therefore, using Sweedler's sumless notation,
$\Delta_q^{(n-1)}(a\varphi^q)=a_{(0)}\varphi^q\otimes\dots\otimes a_{(n-1)}\varphi^q$. Hence to prove that the sections  $q\mapsto\sigma^q(a\varphi^q)\in A_q$ are continuous for all $a\in\C[G]$ it suffices to show that the sections $q\mapsto\sigma^q_i(a\varphi^q)$ of the field $(C(SU_{q^{d_i}}(2)))_{q>0}$ are continuous for all $a\in\C[G]$ and $1\le i\le r$.

Fix $i$ and a continuous family of $*$-isomorphisms $\theta^q\colon \U(SU_q(2))\to\U(SU(2))$ such that $\theta^1=\iota$.
To simplify the notation assume $d_i=1$, so $\sigma^q_i$ is defined by an embedding $\rho^q_i\colon U_q\sltwo\to U_q\g$. To finish the proof it suffices to show that there exists a continuous family of $*$-automorphisms $\gamma^q$ of $\U(G)$ such that $\gamma^q\varphi^q\rho^q_i=\rho_i\theta^q$. Indeed, then exactly as in the first part of the proof, a section $q\mapsto\sigma^q_i(a \varphi^q)$ is continuous if and only if $q\mapsto\sigma^q_i(a\gamma^q\varphi^q)$ is continuous. Since
$$
\sigma^q_i(a\gamma^q\varphi^q)=a\gamma^q\varphi^q\rho^q_i=a\rho_i\theta^q,
$$
the latter section is indeed continuous by definition of the continuous field structure on $(C(SU_q(2)))_{q>0}$.

The automorphisms $\gamma^q$ will be defined by a family of automorphisms $\gamma^q_\lambda$ of $B(V_\lambda)$. Fix $\lambda\in P_+$.
Let $N\in\N$ be such that $\omega(h_i)\le N$ for every $\omega\in P$ such that $V_\lambda(\omega)\ne0$. For every $q>0$ consider the direct sum $\oplus_{s\le N/2}V^q_s$ of $U_q\sltwo$-modules and the corresponding surjective homomorphism $\alpha^q$ from~$U_q\sltwo$ into the algebra
$B_q=\oplus_{s\le N/2}B(V^q_s)$. Then the representation $\pi^q_\lambda\rho^q_i\colon U_q\sltwo\to B(V^q_\lambda)$ factors through $B_q$, so $\pi^q_\lambda\rho^q_i=\beta^q\alpha^q$ for a unique $\beta^q\colon B_q\to B(V_\lambda^q)$. We summarize all the maps involved in the following diagram, which is commutative along solid lines:
\begin{equation} \label{ediag}
\xymatrix{\U(G_q)\ar[r]^{\pi^q_\lambda\ \ } & B(V^q_\lambda)\ar@{-->}[r]^{\varphi^q_\lambda} & B(V_\lambda) & \U(G)\ar[l]_{\ \ \ \pi_\lambda}\\
 & B_q\ar[r]^{\tilde\theta^q}\ar[u]^{\beta^q} & B\ar[u]_{\beta} & \\
\U(SU_q(2))\ar[uu]^{\rho^q_i}\ar[rrr]^{\theta^q}\ar[ur]^{\alpha^q} & & & \U(SU(2))\ar[uu]_{\rho_i}\ar[ul]_\alpha
},
\end{equation}
where $\varphi^q_\lambda\colon B(V^q_\lambda)\to B(V_\lambda)$ and $\tilde\theta^q\colon B_q\to B$ are the isomorphisms defined by the isomorphisms $\varphi^q\colon \U(G_q)\to \U(G)$ and $\theta^q\colon \U(SU_q(2))\to \U(SU(2))$, respectively. Consider the family of homomorphisms $\varphi^q_\lambda\beta^q(\tilde\theta^q)^{-1}\colon B\to B(V_\lambda)$. Since the families of homomorphisms $\tilde\theta^q\alpha^q=\alpha\theta^q\colon U_q\sltwo\to B$ and $\varphi^q_\lambda\beta^q\alpha^q=\varphi^q_\lambda\pi^q_\lambda\rho^q_i=\pi_\lambda\varphi^q\rho^q_i\colon U_q\sltwo\to B(V_\lambda)$ are continuous in the sense defined earlier, and the homomorphisms in the first family are surjective, it follows that the homomorphisms $\varphi^q_\lambda\beta^q(\tilde\theta^q)^{-1}\colon B\to B(V_\lambda)$ depend continuously on $q$. Furthermore, for $q=1$ we get the homomorphism $\beta$. Hence, by a standard result on homomorphisms of finite dimensional C$^*$-algebras, we can choose a continuous family of $*$-automorphisms $\gamma^q_\lambda$ of $B(V_\lambda)$ such that $\gamma^q_\lambda\varphi^q_\lambda\beta^q(\tilde\theta^q)^{-1}
=\beta$ for all $q>0$. In other words, if we replace $\varphi^q_\lambda$ by $\gamma^q_\lambda\varphi^q_\lambda$ in~\eqref{ediag}, we get a commutative diagram, which is what we need.
\ep

\begin{remark} \label{rdense}
The space of sections of the form $q\mapsto a\varphi^q$, $a\in \C[G]$, is not an algebra and is not closed under involution.
But the space of finite sums of sections of the form $q\mapsto f(q)a\varphi^q$, where $a\in\C[G]$ and $f$ is a continuous function, is a $*$-algebra. Indeed, assume $a$ is a matrix coefficient of a finite dimensional representation of $G$ and $\{a_i\}_i$ is a basis in the space spanned by the matrix coefficients of the contragradient representation. Then $(a\varphi^q)^*=\sum_i f_i(q)a_i\varphi^q$ for uniquely defined functions $f_i$. Since $q\mapsto(a\varphi^q)^*$ is a continuous section, the functions $f_i$ must be continuous. Therefore the space is closed under involution. Similarly we check that the space is closed under multiplication. This of course can also be checked without relying on the above theorem. Note also that this space does not depend on the choice of $\varphi^q$.
\end{remark}

For $b>a>0$ denote by $C(G_{[a,b]})$ the C$^*$-algebra of continuous sections of the field $(C(G_q))_{q\in[a,b]}$.

\begin{theorem}[\cite{Na}] \label{tNagy}
For any $b>a>0$ and $q\in[a,b]$ the evaluation map $ev_q\colon C(G_{[a,b]})\to C(G_q)$ is a KK-equivalence.
\end{theorem}

Since this is not exactly how the result is formulated in~\cite[Corollaries~3.8 and~3.11]{Na}, some comments are in order. First, it is stated only for $G=SU(N)$ and it is mentioned that it can be similarly proved for other classical simple compact Lie groups. Once the family of C$^*$-algebras $C(G_q)$ is given a continuous field structure as described above, the general case is essentially identical to $G=SU(N)$. Second, Nagy works with E-theory rather than with KK-theory. Since all the algebras involved are nuclear, there is no difference. Finally, the result in~\cite{Na} is stated only for intervals of the form $[a,1]$ and the evaluations at the end points. In other words, the C$^*$-algebras $I_{a,b}$ (resp. $J_{a,b}$) of continuous sections of $(C(G_q))_{q\in[a,b]}$ vanishing at $a$ (resp., at $b$) are KK-contractible for $b=1$ and all $0<a<1$. From the proof of KK-contractibility of $I_{a,1}$ (which is easier than that of $J_{a,1}$, see \cite[Corollary~3.8]{Na}), it also follows that the algebras $I_{a,b}$ and $J_{a,b}$ are KK-contractible for $b<1$. Using the canonical isomorphisms $G_q\cong G_{q^{-1}}$ we therefore conclude that $I_{a,b}$ and $J_{a,b}$ are KK-contractible for $a<b\le1$ and $1\le a <b$. For $a<1<b$ the exact sequences $0\to I_{1,b}\to I_{a,b}\to I_{a,1}\to0$ and $0\to J_{a,1}\to J_{a,b}\to J_{1,b}\to0$ show that $I_{a,b}$ and $J_{a,b}$ are KK-contractible for all $b>a>0$. Therefore $ev_q\colon C(G_{[a,b]})\to C(G_q)$ is a KK-equivalence for $q=a,b$. Finally, for $q\in(a,b)$ the kernel of $ev_q$ is $J_{[a,q]}\oplus I_{[q,b]}$, so $ev_q$ is again a KK-equivalence.

\bigskip

\section{Continuous family of Drinfeld twists}
\label{stwist}

As in the previous section, fix a continuous family of $*$-isomorphisms $\varphi^q\colon \U(G_q)\to\U(G)$ with $\varphi^1=\iota$.

For $q>0$, let $\hbar_q\in i\R$ be such that $q=e^{\pi i\hbar_q}$. Denote by $t\in\g\otimes\g$ the
$\g$-invariant element defined by the $\ad$-invariant symmetric form on $\g$ such that the induced form on $\h^*$ satisfies $(\alpha_i, \alpha_j)=d_ia_{ij}$. By a result of Kazhdan and Lusztig~\cite{KL1}, see~\cite{NT3} for details, for every $q>0$ there exists  a unitary element $\F^q\in \U(G\times G)$, which we call a unitary Drinfeld twist, such that
 \enu{i}
$(\varphi^q\otimes\varphi^q)\Dhat_q=\F^q\Dhat\varphi^q(\cdot){\F^q}^*$; \enu{ii}
$(\hat\eps\otimes\iota)(\F^q)=(\iota\otimes\hat\eps)(\F^q)=1$; \enu{iii}
$(\varphi^q\otimes\varphi^q)(\RR_q)=\F^q_{21}q^t{\F^q}^*$, where $\RR_q\in \U(G_q\times G_q)$ is the universal $R$-matrix; \enu{iv}
$(\iota\otimes\Dhat)({\F^q})^*
(1\otimes{\F^q})^*(\F^q\otimes1)(\Dhat\otimes\iota)(\F^q)=\Phi_{KZ}(\hbar_q t_{12},\hbar_q t_{23})$, where $\Phi_{KZ}$ is Drinfeld's KZ-associator.

Such an element is not unique, but by~\cite{NT4} any other unitary Drinfeld twist (for the same isomorphism $\varphi^q$) has the form $(c\otimes c)\F^q\Dhat(c)^*$ for a unitary element $c$ in the center of $\U(G)$.

We say that a family $\{\F^q\}_{q>0}$ of unitary Drinfeld twists is continuous if the map $$q\mapsto\F^q\in W^*(G\times G)$$ is continuous in the strong operator topology on the von Neumann algebra $W^*(G\times G)\subset \U(G\times G)$ of $G\times G$. In other words, the map $q\mapsto(\pi_\lambda\otimes\pi_\nu)(\F^q)$ is continuous for all $\lambda,\nu\in P_+$.

\begin{theorem}
There exists a continuous family of unitary Drinfeld twists $\F^q$ such that $\F^1=1$. Furthermore, if $\{\psi^q\colon \U(G_q)\to\U(G)\}_{q>0}$ is another continuous family of $*$-isomorphisms such that $\psi^1=\iota$, and $\{\E^q\}_{q>0}$ is a corresponding continuous family of unitary Drinfeld twists with $\E^1=1$, then there exists a unique continuous family of unitary elements $u^q\in \U(G)$ such that
$$
u^1=1,\ \ \hbox{and} \ \ \psi^q=u^q\varphi^q(\cdot){u^q}^*\ \ \hbox{and}\ \ \E^q=(u^q\otimes u^q)\F^q\Dhat(u^q)^*\ \ \hbox{for all} \ \ q>0.
$$
\end{theorem}

\bp To prove existence, consider the set $\Omega$ of pairs $(q,\F)$, where $q>0$ and $\F$ is a unitary Drinfeld twist for $\varphi^q$. It is a closed subset of the direct product of $\R^*_+$ and the unitary group of $W^*(G\times G)$ (this is used already in the proof of~\cite[Lemma~3.2]{NT3}), so it is a locally compact space. Let $p\colon \Omega\to \R^*_+$ be the projection onto the first coordinate. The compact abelian group of elements of the form $(c\otimes c)\Dhat(c)^*$, where $c$ is a unitary element in the center of $\U(G)$, acts freely by multiplication on the right on $\Omega$, and by~\cite[Theorem~5.2]{NT4} this action is transitive on each fiber of the map $p$. Therefore if this group were a compact Lie group, then by a theorem of Gleason~\cite{Gl}, $p\colon \Omega\to \R^*_+$ would be a fiber bundle, hence $p$ would have a continuous section. Since the group of elements of the form $(c\otimes c)\Dhat(c)^*$ is not a Lie group, we cannot apply Gleason's theorem directly and will proceed as follows.

Choose an increasing sequence of finite subsets $P_n\subset P_+$ such that $P_1=\{0\}$ and $\cup_n P_n=P_+$. For every $q>0$ we will construct a sequence of unitary Drinfeld twists $\F^q_n$ such that $\F_n^1=1$, the map $q\mapsto(\pi_\lambda\otimes\pi_\nu)(\F^q_n)$ is continuous for all $\lambda,\nu\in P_n$ and $n\ge1$, and $(\pi_\lambda\otimes\pi_\nu)(\F^q_{n+1})=(\pi_\lambda\otimes\pi_\nu)(\F^q_n)$ for all $\lambda,\nu\in P_n$ and $n\ge1$. Then, for every $q>0$, the sequence $\{\F^q_n\}_n$ converges to a unitary Drinfeld twist $\F^q$ with the required properties.

For $n=1$ and $q\ne1$ we take $\F^q_1$ to be any unitary Drinfeld twist, and we take $\F^1_1=1$.

Assume the Drinfeld twists $\F^q_n$ are already constructed for some $n$. Denote by $\Omega_{n+1}$ the set of pairs $(q, U)$, where $U=(U_{\lambda,\nu})_{\lambda,\nu}$ is a unitary element in~$\prod_{(\lambda,\nu)\in P_{n+1}\times P_{n+1}\setminus P_n\times P_n}B(V_\lambda\otimes V_\nu)$ such that there exists a unitary Drinfeld twist $\F$ for $\varphi^q$ satisfying
\begin{align*}
(\pi_\lambda\otimes\pi_\nu)(\F)&=U_{\lambda,\nu}\ \ \hbox{for all}\ \ (\lambda,\nu)\in P_{n+1}\times P_{n+1}\setminus P_n\times P_n, \\
(\pi_\lambda\otimes\pi_\nu)(\F)&=(\pi_\lambda\otimes\pi_\nu)(\F^q_n) \ \ \hbox{for all}\ \ \lambda,\nu\in P_n.
\end{align*}
Let $p_{n+1}\colon\Omega_{n+1}\to\R^*_+$ be the projection onto the first coordinate. The set $\Omega_{n+1}$ is a closed subset of the direct product of $\R^*_+$ and the unitary group of $\prod_{(\lambda,\nu)\in P_{n+1}\times P_{n+1}\setminus P_n\times P_n}B(V_\lambda\otimes V_\nu)$. For every $q>0$ the fiber $p_{n+1}^{-1}(q)$ is nonempty, as it contains the element $((\pi_\lambda\otimes\pi_\nu)(\F^q_n))_{\lambda,\nu}$.

Let $S_{n+1}$ be the set of weights $\lambda\in P_+$ such that either $\lambda\in P_{n+1}$, or $V_\lambda$ is equivalent to a subrepresentation of $V_\nu\otimes V_\eta$ for some $\nu,\eta\in P_{n+1}$, in which case we write $V_\lambda\prec V_\nu\otimes V_\eta$. Let $K_{n+1}=\prod_{\lambda\in S_{n+1}}\T$. We have a homomorphism $\rho_{n+1}$ from $K_{n+1}$ into the unitary group of
$\prod_{(\lambda,\nu)\in P_{n+1}\times P_{n+1}\setminus P_n\times P_n}B(V_\lambda\otimes V_\nu)$: $\rho_{n+1}(c)$ acts on the isotypic component of $V_\nu\otimes V_\eta$ of type $V_\lambda$ as multiplication by $c_\nu c_\eta \bar c_\lambda$. We also have a similar homomorphism $\theta_{n+1}$ from $K_{n+1}$ into the unitary group of $\prod_{\lambda,\nu\in P_n}B(V_\lambda\otimes V_\nu)$.

The group $\ker \theta_{n+1}$ acts on $\Omega_{n+1}$ by multiplication by $\rho_{n+1}(c)$ on the right. On every fiber of $p_{n+1}$ this action is transitive, and the stabilizer of every point is $\ker\rho_{n+1}\cap\ker\theta_{n+1}$. Since $\ker \theta_{n+1}$ is a compact Lie group, by Gleason's theorem we conclude that $p_{n+1}\colon\Omega_{n+1}\to \R^*_+$ is a fiber bundle, hence it is a trivial bundle. Choosing a continuous section of this bundle, by definition of $\Omega_{n+1}$ we conclude that there exist unitary Drinfeld twists $\E^q$ such that the map  $q\mapsto(\pi_\lambda\otimes\pi_\nu)(\E^q)$ is continuous for all $\lambda,\nu\in P_{n+1}$ and  $(\pi_\lambda\otimes\pi_\nu)(\E^q)=(\pi_\lambda\otimes\pi_\nu)(\F^q_n)$ for all $\lambda,\nu\in P_n$. There exists a unitary central element $c$ in~$\U(G)$ such that $\E^1=(c^*\otimes c^*)\Dhat(c)$. We can then set $\F^q_{n+1}=\E^q(c\otimes c)\Dhat(c)^*$. This finishes the proof of existence.

\smallskip

Assume now that $\{\psi^q\colon \U(G_q)\to\U(G)\}_{q>0}$ is another continuous family of $*$-isomorphisms such that $\psi^1=\iota$, and $\{\E^q\}_{q>0}$ is a corresponding continuous family of unitary Drinfeld twists with $\E^1=1$. For every $\lambda\in P_+$, let $\varphi^q_\lambda,\psi^q_\lambda\colon B(V^q_{\lambda})\to B(V_\lambda)$ be the isomorphisms defined by $\varphi^q$ and $\psi^q$. The set of unitaries $v\in B(V_\lambda)$ such that $\psi^q_\lambda=v\varphi^q_\lambda(\cdot)v^*$ forms a circle bundle over $\R^*_+$, so it has a continuous section $v^q_\lambda$. Since $\psi^1_\lambda=\varphi^1_\lambda=\iota$, we may assume that $v^1_\lambda=1$. The unitaries $v^q_\lambda$ define a continuous family of unitaries $v^q\in \U(G)$.

For every $q>0$, the element $(v^q\otimes v^q)\F^q\Dhat(v^q)^*$ is a unitary Drinfeld twist for $\psi^q$. Hence, for every $q$, there exists a unitary central element $c\in \U(G)$ such that
\begin{equation} \label{ecohom}
\E^q=(v^q\otimes v^q)\F^q\Dhat(v^q)^*(c\otimes c)\Dhat(c)^*.
\end{equation}
Furthermore, the element $c$ is defined up to a group-like unitary element in the center of $\U(G)$, that is, up to an element of the center $Z(G)$ of $G$. Therefore, applying once again Gleason's theorem (which in this case is quite obvious as $Z(G)$ is finite), we see that the set of pairs $(q,c)$ with $c$ satisfying \eqref{ecohom} is a principle $Z(G)$-bundle over $\R^*_+$, hence it has a continuous section $q\mapsto (q,c^q)$. The element $c^1$ is group-like, so replacing $c^q$ by $c^q{c^1}^*$ we may assume that $c^1=1$. Letting $u^q=c^qv^q$, we get the required continuous family of unitary elements.

Finally, if $\tilde u^q$ is another continuous family of unitary elements with the same properties, then $c^q=\tilde u^q{u^q}^*$ is a unitary central group-like element in~$\U(G)$, hence $c^q\in Z(G)$. Since $c^q$ depends continuously on $q$, $Z(G)$ is finite and $c^1=1$, we conclude that $c^q=1$ for all $q$.
\ep

Another way of formulating the above result is to say that the set of triples $(q,\varphi,\F)$ such that $q>0$, $\varphi\colon\U(G_q)\to\U(G)$ is a $*$-isomorphism extending the canonical identification of the centers and~$\F$ is a unitary Drinfeld twist for $\varphi$, has a structure of a principal $U(W^*(G))/Z(G)$-bundle over~$\R^*_+$, where $U(W^*(G))$ is the unitary group of $W^*(G)$.

\smallskip

Note also that by analyzing the proof of Kazhdan and Lusztig~\cite{KL2} one could hope to prove a stronger result: the family of unitary Drinfeld twists $\F^q$ can be chosen to be real-analytic for an appropriate choice of $\varphi^q$.

\bigskip

\section{Family of Dirac operators}\label{sdirac}

We continue by fixing a continuous family of $*$-isomorphisms $\varphi^q\colon \U(G_q)\to\U(G)$ with $\varphi^1=\iota$
and a continuous family of unitary Drinfeld twists $\F^q$ for $\varphi^q$ with $\F^1=1$.

For every $q>0$ we have a Dirac operator $D_q$ on $G_q$ defined as follows~\cite{NT2}.

Consider a basis~$\{x_i\}_i$ of $\g$ such
that $(x_i ,x_j)=-\delta_{ij}$, and let $\gamma\colon\g\to\Clg$ denote
the inclusion of $\g$ into the complex Clifford algebra with the
convention that $\gamma (x_i)^2=-1$. Identifying ${\mathfrak{so} }(\g)$
with ${\mathfrak{spin}}(\g)$, the adjoint action is defined by the
representation $\add\colon\g\to{\mathfrak{spin}} (\g)\subset\Clg$ given
by
$$\add(x)=\frac{1}{4}\sum_i \gamma (x_i )\gamma ([x,x_i ]).$$
We denote by the same symbol $\add$ the corresponding homomorphism
$\U(G)\to\Clg$.

Let $s\colon\Clg\to \End(\Sp)$ be an irreducible representation. Denote
by $\partial$ the representation of $U\g$ by left-invariant
differential operators. Identifying the sections $\Gamma (S)$ of the
spin bundle $S$ over $G$ with $C^\infty (G)\otimes\Sp$, the Dirac
operator $D\colon C^\infty (G)\otimes\Sp\to C^\infty (G)\otimes\Sp$
defined using the Levi-Civita connection, can be written as
$D=(\partial\otimes s)(\D)$, where $\D\in U\g\otimes\Clg$ is given by
the formula
$$\D=\sum_i (x_i\otimes\gamma (x_i )+\frac{1}{2}\otimes\gamma(x_i)\add(x_i)).$$


\smallskip

Let $(L^2 (G_q),\pi_{q,r},\xi^q_h)$ be the GNS-triple defined by
the Haar state on $C(G_q)$. The right regular
representation of the von Neumann algebra~$W^*(G_q)\subset\U(G_q)$ of $G_q$ on $L^2(G_q)$, denoted by $\partial_q$, is defined by

\begin{equation*} \label{eRegL}
\partial_q(\omega)\pi_{r,q}(a)\xi^q_h
=(\pi_{r,q}\otimes\omega)\Delta_q(a)\xi^q_h
=a_{(1)}(\omega)\pi_{r,q}(a_{(0)})\xi^q_h.
\end{equation*}

The Dirac operator $D_q$ on $G_q$ is the unbounded operator on
$L^2(G_q)\otimes\Sp$ defined by
$$
D_q=(\partial_q\otimes s)(\D_q),
$$
where
$\D_q\in\U(G_q)\otimes\Clg$ is given by
$$
\D_q=(\varphi^q\otimes\iota)^{-1}((\iota\otimes\add)(\F^q)
\D(\iota\otimes\add)({\F^q})^*).
$$

Our goal is to show that the family $(D_q)_q$ is continuous in the sense that it defines a Kasparov $(C(G_{[a,b]}),C[a,b])$-module.

\begin{lemma}
The family $(L^2(G_q))_{q>0}$ has a unique structure of a continuous field of Hilbert spaces such that the vector field $q\mapsto \pi_{r,q}(a^q)\xi^q_h$ is continuous for every continuous section $q\mapsto a^q$ of the field $(C(G_q))_{q>0}$.
\end{lemma}

\bp It suffices to show that the function $q\mapsto (\pi_{r,q}(a^q)\xi^q_h,\xi^q_h)$ is continuous. But this is clear, since if $a^q=a\varphi^q$, where $a$ is a matrix coefficient of a nontrivial irreducible representation of $G$, then the function is zero by the orthogonality relations.
\ep

It is easy to see that the continuous field $(L^2(G_q))_{q>0}$ is trivial. To formulate a more precise result, recall the exact form of the orthogonality relations.

First let us introduce some notation. For a weight $\beta=\sum_i c_i\alpha_i$ ($c_i\in\C$) put $h^q_\beta=\sum_i c_id_ih_i^q\in\U(G_q)$ and $K^q_\beta=q^{h^q_\beta}$; note that for $q=1$ the element $h_\beta$ is characterized by $\lambda(h_\beta)=(\lambda,\beta)$ for any weight~$\lambda$, and for $q\ne1$ we have $K^q_{\alpha_i}=K^q_i$. Let $\rho$ be half the sum of positive roots. Then $K^q_{2\rho}$ is the Woronowicz character $f_{-1}$ for $G_q$; in particular, for the square of the antipode $\hat S_q$ on $\U(G_q)$ we have $\hat S_q^2(\omega)=K^q_{-2\rho}\omega K^q_{2\rho}$. Put $\dim_q(V^q_\lambda)=\Tr(\pi^q_\lambda(K^q_{2\rho}))=\Tr(\pi^q_\lambda(K^q_{-2\rho}))$.

For $\xi,\zeta\in V^q_\lambda$ define $a^{\lambda;q}_{\xi,\zeta}\in \C[G_q]$ by
$$
a^{\lambda;q}_{\xi,\zeta}(\omega)=(\pi^q_\lambda(\omega)\zeta,\xi)\ \ \hbox{for}\ \ \omega\in\U(G_q).
$$
Then the orthogonality relations state that the vectors $\pi_{r,q}(a^{\lambda;q}_{\xi,\zeta})\xi^q_h$ are mutually orthogonal for different $\lambda$, and
\begin{equation*} \label{eorth}
(\pi_{r,q}(a^{\lambda;q}_{\xi,\zeta})\xi^q_h,\pi_{r,q}(a^{\lambda;q}_{\xi',\zeta'})\xi^q_h)=\frac{1}{\dim_q(V^q_\lambda)}
\overline{(\pi^q_\lambda(K^q_{2\rho})\xi,\xi')}(\zeta,\zeta').
\end{equation*}

Let $d^q\in\U(G_q)$ be the element such that
$$
\pi^q_\lambda(d^q)=\frac{\dim(V_\lambda)^{1/2}}{\dim_q(V^q_\lambda)^{1/2}}\pi^q_\lambda(K^q_\rho)\ \ \hbox{for all}\ \ \lambda\in P_+.
$$

\begin{lemma} \label{lintert}
The linear operator $W_q\colon \pi_{r,q}(\C[G_q])\xi^q_h\to  \pi_{r}(\C[G])\xi_h$ defined by
$$
W_q\pi_{r,q}(a\varphi^q)\xi_h^q=a_{(0)}(\varphi^q(d^q))\pi_r(a_{(1)})\xi_h\ \ \hbox{for}\ \ a\in\C[G]
$$
is unitary. It has the property
$$
W_q\partial_q(\omega)=\partial(\varphi^q(\omega)) W_q \ \ \hbox{for all}\ \ \omega\in\U(G_q).
$$
\end{lemma}

\bp By the orthogonality relations we have a decomposition
$L^2(G_q)=\oplus_{\lambda\in P_+}\overline{V^q_\lambda}\otimes V^q_\lambda$, with $\pi_{r,q}(a^{\lambda;q}_{\xi,\zeta})\xi^q_h\in L^2(G_q)$ corresponding to $ \dim_q(V^q_\lambda)^{-1/2}\overline{\pi^q_\lambda(K^q_{\rho})\xi}\otimes\zeta\in \overline{V^q_\lambda}\otimes V^q_\lambda$. The isomorphism $\varphi^q\colon\U(G_q)\to\U(G)$ is implemented by unitaries $u^q_\lambda\colon V^q_\lambda\to V_\lambda$. Then the unitaries
$$
\overline{u^q_\lambda}\otimes u^q_\lambda\colon \overline{V^q_\lambda}\otimes V^q_\lambda\to \overline{V_\lambda}\otimes V_\lambda
$$
define a unitary $L^2(G_q)\to L^2(G)$. This is exactly the unitary $W_q$, since if $a=a^\lambda_{\xi,\zeta}\in\C[G]$ then $a\varphi^q=a^{\lambda;q}_{{u^q_\lambda}^*\xi,{u^q_\lambda}^*\zeta}$ and hence the vector $\pi_{r,q}(a\varphi^q)\xi_h^q$ is mapped onto
$$
\frac{\dim(V_\lambda)^{1/2}}{\dim_q(V^q_\lambda)^{1/2}}\pi_{r}(a^{\lambda}_{u^q_\lambda\pi^q_\lambda(K^q_\rho){u^q_\lambda}^*\xi,\zeta})\xi_h,
$$
which gives the formula in the formulation, as
$$
a^{\lambda}_{u^q_\lambda\pi^q_\lambda(K^q_\rho){u^q_\lambda}^*\xi,\zeta}=(\pi_\lambda(\cdot)\zeta,\pi_\lambda(\varphi^q(K^q_\rho))\xi)=
a(\varphi^q(K^q_\rho)\,\cdot)=a_{(0)}(\varphi^q(K^q_\rho))a_{(1)}.
$$

The last statement in the formulation follows either by a direct computation or by observing that~$\partial_q(\omega)$ acts on $\overline{V^q_\lambda}\otimes V^q_\lambda\subset L^2(G_q)$ as $1\otimes\pi^q_\lambda(\omega)$.
\ep

The extension of $W_q$ to $L^2(G_q)$ we continue to denote by the same symbol. Therefore the unitaries~$W_q$ define an isomorphism of the continuous field $(L^2(G_q))_{q>0}$ onto the trivial field with fiber~$L^2(G)$.

\smallskip

For $b>a>0$ denote by $L^2(G_{[a,b]})$ the right Hilbert $C[a,b]$-module of continuous sections of $(L^2(G_q))_{q\in[a,b]}$. The C$^*$-algebra $C(G_{[a,b]})$ acts on $L^2(G_{[a,b]})$ via GNS-representations $\pi_{r,q}$; we denote by $\pi$ the corresponding homomorphism from $C(G_{[a,b]})$ into the algebra of adjointable operators on~$L^2(G_{[a,b]})$. The operators $D_q$ define an unbounded operator $D_{[a,b]}$ on $L^2(G_{[a,b]})\otimes \Sp$ with domain of definition consisting of continuous vector fields $\xi$ such that $\xi(q)\in\Dom(D_q)$ for all $q$ and the vector field $q\mapsto D_q\xi(q)$ is continuous. When $G$ is even-dimensional we also define a grading on $L^2(G_{[a,b]})\otimes \Sp$ using the chirality element $\chi\in\Clg$.

\begin{theorem}
For any $b>a>0$ the triple $(L^2(G_{[a,b]})\otimes\Sp,\pi(\cdot)\otimes1,D_{[a,b]})$ is an unbounded Kasparov $(C(G_{[a,b]}),C[a,b])$-module of the same parity as the dimension of $G$.
\end{theorem}

\bp By definition of an unbounded Kasparov module \cite{BJ} it suffices to check that
\begin{itemize}
\item[(i)] $D_{[a,b]}$ is a regular self-adjoint operator such that $(1+D^2_{[a,b]})^{-1}$ is generalized compact (recall that regularity means that the operator $1+D^2_{[a,b]}$ is surjective);
\item[(ii)] there exists a dense $*$-subalgebra $\A$ of $C(G_{[a,b]})$ and an $\A$-invariant core of $D_{[a,b]}$ such that the commutators $[D_{[a,b]},\pi(c)\otimes1]$ are bounded on this core for all $c\in\A$.
\end{itemize}

By definition of $D_q$ and Lemma \ref{lintert} we have
\begin{equation} \label{eequiv}
(\partial\otimes s\,\add)(\F^q)^*(W_q\otimes1)D_q=D(\partial\otimes s\,\add)(\F^q)^*(W_q\otimes1).
\end{equation}
Therefore the unitaries $(\partial\otimes s\,\add)(\F^q)^*(W_q\otimes1)$ define an isomorphism of the continuous field $(L^2(G_q)\otimes\Sp)_q$ onto the constant field with fiber $L^2(G)\otimes\Sp$, which maps $D_{[a,b]}$ onto the operator which acts as $D$ on every fiber. This immediately gives (i).

\smallskip

To prove (ii) consider the space $\A\subset C(G_{[a,b]})$ of finite sums of sections of the form $q\mapsto f(q)c\,\varphi^q$ with $f$ a continuous function and $c\in\C[G]$. By Remark \ref{rdense} this is a $*$-algebra. It is dense in~$C(G_{[a,b]})$ since it is dense in every fiber and is closed under multiplication by continuous functions. The linear span of vectors fields of the form $q\mapsto\pi_{r,q}(c^q)\xi^q_h\otimes\zeta$, where $c\in\A$ and $\zeta\in\Sp$, is an $\A$-invariant core for $D_{[a,b]}$. This follows e.g.~from \eqref{eequiv} and the fact that vectors of the form $\pi_r(c)\xi^h\otimes\zeta$, where  $c\in\C[G]$ and $\zeta\in\Sp$, span a core for $D$ which is a union of $D$-invariant finite dimensional subspaces.

It remains to check boundedness of commutators. For this it is enough to show that for every $c\in\C[G]$ the commutators $[D_q,\pi_{r,q}(c\,\varphi^q)\otimes1]$ are uniformly bounded on $[a,b]$. By \cite[Proposition~3.1]{NT2} we have
$$
[D_q,\pi_{r,q}(c\,\varphi^q)\otimes1]=-(\pi_{r,q}(c_{(0)}\varphi^q)\otimes1)
(\partial_q\circ(\varphi^q)^{-1}\otimes s)
(c_{(1)}\otimes\iota\otimes\iota)(U_qT_qU^*_q),
$$
where
$$
U_q=(\iota\otimes\iota\otimes\add)((\F^q\otimes1)
(\Dhat\otimes\iota)(\F^q))\in W^*(G)\bar\otimes W^*(G)\otimes\Clg
$$
and $T_q\in\U(G\times G)\otimes\Clg$ is defined by
$$
T_q=(\iota\otimes\iota\otimes\gamma)(t_{13})
+(\iota\otimes\iota\otimes\gamma)(t_{23})
-(\iota\otimes\iota\otimes\add)({\Phi^q})^*
(\iota\otimes\iota\otimes\gamma)(t_{23})
(\iota\otimes\iota\otimes\add)(\Phi^q),
$$
with $\Phi^q=\Phi_{KZ}(\hbar_q t_{12},\hbar_q t_{23})$. It is therefore enough to prove that the operators
$$
(c\otimes\iota\otimes\iota)(U_qT_qU^*_q)\in\U(G)\otimes\Clg
$$
are uniformly bounded on $[a,b]$ for all $c\in\C[G]$. Equivalently, the operators $(\pi\otimes\iota\otimes\iota)(T_q)$ are uniformly bounded for any finite dimensional unitary representation $\pi$ of $G$, which, it turn, is the same as uniform boundedness of
$$
[(\pi\otimes\iota\otimes\gamma)(t_{23}),(\pi\otimes\iota\otimes\add)(\Phi_{KZ}(\hbar_q t_{12},\hbar_q t_{23}))].
$$
By \cite[Proposition~3.6]{NT2} the latter property indeed holds: the norm of the above commutator is bounded by~$6\|(\pi\otimes\gamma)(t)\|$ independently of $q>0$.
\ep

Recall that by Theorem \ref{tNagy} the evaluation map $ev_q\colon C(G_{[a,b]})\to C(G_q)$ is a KK-equivalence for any $q\in[a,b]$. For $q,q'\in[a,b]$ define an invertible element $\gamma_{q,q'}$ in~$KK(C(G_q),C(G_{q'}))$ by $\gamma_{q,q'}=[ev_q]^{-1}[ev_{q'}]$. It does not depend on the segment $[a,b]$ containing $q$ and $q'$. Denote also by~$[D_q]$ the element of $KK_i(C(G_q),\C)$ defined by $D_q$, where $i=\dim G \mod 2$.

\begin{corollary}
For any $q,q'>0$ we have $\gamma_{q,q'}[D_{q'}]=[D_q]$.
\end{corollary}

\bp
Denote by $[D_{[a,b]}]\in KK_i(C(G_{[a,b]}),C[a,b])$ the class of the Kasparov module $$(L^2(G_{[a,b]})\otimes\Sp,\pi(\cdot)\otimes1,D_{[a,b]}).$$ If $\tilde{ev}_q\colon C[a,b]\to\C$ is the evaluation at $q$ then clearly $[D_{[a,b]}]\,[\tilde{ev}_q]=[ev_q]\,[D_q]$. Since $[\tilde{ev}_q]$ does not depend on $q$, we thus see that the class $[ev_q]\,[D_q]\in KK_i(C(G_{[a,b]}),\C)$ does not depend on $q\in[a,b]$ either, which is what we need.
\ep

\end{document}